\documentclass[a4paper,12pt]{article}
\usepackage{amsmath}
\usepackage{amsthm}
\usepackage{amssymb}
\usepackage{amscd}

\title{Spherical harmonic polynomials for higher bundles}
\author{Yasushi Homma\thanks{Department of Mathematical Sciences, Waseda University, 
  3-4-1 Ohkubo, Shinjuku-ku, Tokyo, 169-8555, JAPAN. \endgraf 
  \textit{e-mail address}: homma@gm.math.waseda.ac.jp}}
\date{}

\theoremstyle{plain}
\newtheorem{thm}{Theorem}[section]
\newtheorem{prop}[thm]{Proposition}
\newtheorem{cor}[thm]{Corollary}
\newtheorem{lem}[thm]{Lemma}
\theoremstyle{definition}

\numberwithin{equation}{section}

\theoremstyle{remark}
\newtheorem{rem}{Remark}[section]
\newtheorem{exam}{Example}[section]

\newcommand{\id}{\mathrm{id}}
\newcommand{\Ad}{\mathrm{Ad}}

\allowdisplaybreaks[3]
\begin{document}
\maketitle
\begin{abstract}
We give a method of decomposing bundle-valued polynomials compatible with the action of the Lie group $Spin(n)$, where important tools are $Spin(n)$-equivariant operators and their spectral decompositions. In particular, the top irreducible component is realized as an intersection of kernels of these operators. 
\end{abstract}
\section{Introduction}\label{sec:0}
Spherical harmonic polynomials or spherical harmonics are polynomial solutions of the Laplace equation $\square \phi(x)=\sum \partial^2\phi/\partial x_i^2=0$ on $\mathbf{R}^n$. These are fundamental and classical objects in mathematics and physics. It is natural that we consider vector-valued spherical harmonic polynomials. For example, the polynomial solutions of the Dirac equation $D\phi(x)=0$ on $\mathbf{R}^n$ are studied in Clifford analysis (see \cite{DSS}, \cite{H1}, and \cite{T}). They are spinor-valued polynomials and called spherical monogenics. We also have other examples in \cite{Bu}, \cite{Fo}, \cite{H2}, and \cite{IT}, where we can give spectral information of some basic operators on sphere. Recently, the first-order $Spin(n)$-equivariant differential operators have been studied like Dirac operator and Rarita-Schwinger operator (see \cite{Br}-\cite{Bu}, \cite{H3}, and \cite{H4}). These operators are called higher spin Dirac operators or Stein-Weiss operators. In this paper, we give a method to analyze polynomial sections for natural bundles on $\mathbf{R}^n$ by using higher spin Dirac operators and Clifford homomorphisms. Here, Clifford homomorphism is a natural generalization of Clifford algebra given in \cite{H3} and \cite{H4}. 
 
 Let $S^q$ (resp. $H^q$) be the spaces of polynomials (resp. harmonic polynomials) with degree $q$ on the $n$-dimensional Euclidean space $\mathbf{R}^n$. We know that $H^q$ is an irreducible representation space for $Spin(n)$, and $S^q$ has irreducible decomposition, $\oplus_{0\le k \le [q/2]}r^{2k}H^{q-2k}$, where $r$ is $|x|=\sqrt{x_1^2+\cdots+x_n^2}$. To give such a decomposition, we use the invariant operator $-r^2\square$ and its spectral decomposition. In particular, the top component $H^q$ is the kernel of the operator $-r^2\square$. Now, we consider natural irreducible bundle $\mathbf{R}^n\times V_{\rho}$ on $\mathbf{R^n}$, where $V_{\rho}$ is an irreducible representation space with highest weight $\rho$ for $Spin(n)$. Our interest is to analyze the space of $V_{\rho}$-valued polynomials, $S^q\otimes V_{\rho}$. For that purpose, we use higher spin Dirac operators $\{ D^{\rho}_{\lambda_k} \}_k$ and algebraic operators $\{ x^{\rho}_{\lambda_k} \}_k$. Then we have an invariant operator $E$ whose spectral decomposition gives the irreducible decomposition of $S^q\otimes V_{\rho}$ like the operator $-r^2\square$. In particular, the top irreducible component is the kernel of $E$ and realized as an intersection of kernels of higher spin Dirac operators.

\section{Clifford Homomorphisms}\label{sec:1}
In this section, we review Clifford homomorphisms given in \cite{H4}. Let $\mathfrak{spin}(n)\simeq \mathfrak{so}(n)$ be the Lie algebra of the spin group $Spin(n)$ or orthogonal group $SO(n)$. The Lie algebra $\mathfrak{spin}(n)$ is realized by using the Clifford algebra $Cl_n$ associated to $\mathbf{R}^n$: we choose the standard basis $\{e_i\}_i$ of $\mathbf{R}^n$ and put $[e_i,e_j]:=e_ie_j-e_je_i$ in $Cl_n$. Then $\{[e_i,e_j]\}_{i,j}$ span the Lie algebra $\mathfrak{spin}(n)$ in $Cl_n$. 

The irreducible finite dimensional unitary representations of $\mathfrak{spin}(n)$ or $Spin(n)$ are parametrized by dominant weights $\rho=(\rho^1,\cdots,\rho^m) \in \mathbf{Z}^m \cup (\mathbf{Z}+1/2)^m$ satisfying that
\begin{gather}
\rho^1\ge \cdots \ge \rho^{m-1}\ge |\rho^m|, \quad \textrm{for $n=2m$}, \label{eqn:1-6} \\
\rho^1\ge \cdots \ge \rho^{m-1}\ge \rho^m\ge 0, \quad \textrm{for $n=2m+1$}.\label{eqn:1-7}
\end{gather}
We denote by $(\pi_{\rho},V_{\rho})$ not only the representation of $Spin(n)$ but also its infinitesimal one of $\mathfrak{spin}(n)$ with highest weight $\rho$. When writing dominant weights, we denote a string of $j$ $k$'s for $k$ in $\mathbf{Z}\cup (\mathbf{Z}+1/2)$ by $k_j$. For example, the adjoint representation $(\Ad, \mathbf{R}^n\otimes \mathbf{C})$ of $Spin(n)$ (resp. $\mathfrak{spin}(n)$) has the highest weight $(1,0_{m-1})$, where the action is $\pi_{\Ad}(g)u=g u g^{-1}$ for $g$ in $Spin(n)$ (resp. $\pi_{\Ad}([e_i,e_j])u:=[[e_i,e_j], u]$).  

We consider an irreducible representation $(\pi_{\rho},V_{\rho})$ and the tensor representation $(\pi_{\rho}\otimes \pi_{\Ad}, V_{\rho}\otimes_{\mathbf{C}} \mathbf{R}^n)$. We decompose it to irreducible components, $V_{\rho}\otimes_{\mathbf{C}} \mathbf{R}^n=\sum_{0\le k \le N} V_{\lambda_k}$. For $u$ in $\mathbf{R}^n$, we have the following bilinear mapping for each $k$:
\begin{equation}
\mathbf{R}^n\times V_{\rho}\ni (u,\phi)\mapsto p^{\rho}_{\lambda_k}(u)\phi:=\Pi^{\rho}_{\lambda_k}(\phi \otimes u)\in V_{\lambda_k},   \label{eqn:1-1}
\end{equation}
where $\Pi^{\rho}_{\lambda_k}$ is the orthogonal projection from $V_{\rho}\otimes_{\mathbf{C}} \mathbf{R}^n$ onto $V_{\lambda_k}$. We call the linear mapping $p^{\rho}_{\lambda_k}(u):V_{\rho}\to V_{\lambda_k}$ \textit{the Clifford homomorphism} from $V_{\rho}$ to $V_{\lambda_k}$, and denote by $(p^{\rho}_{\lambda_k}(u))^{\ast}$ the adjoint operator of $p^{\rho}_{\lambda_k}(u)$ with respect to the inner products on $V_{\rho}$ and $V_{\lambda_k}$. If we consider the spinor representation $(\pi_{\Delta},V_{\Delta})$, then the Clifford homomorphism from $V_{\Delta}$ to itself is the usual Clifford action of $\mathbf{R}^n$ on $V_{\Delta}$, which satisfy the relation $e_ie_j+e_je_i=-\delta_{ij}$. In general cases, we have a lot of relations among these homomorphisms. 
\begin{thm}[\cite{H4}]\label{relation n}
For any non-negative integer $q$, we define the bilinear mapping $r^q_{\rho}$ as follows: 
\begin{multline}
r^{q}_{\rho}:\mathbf{R}^n\times \mathbf{R}^n \ni (u,v) \mapsto \\
    \left( -\frac{1}{4} \right)^{q} \sum_{l_1,\cdots, l_{q-1}} \pi_{\rho}([u, e_{l_1}])\pi_{\rho}([e_{l_1},e_{l_2}]) \cdots \pi_{\rho}([e_{l_{q-1}}, v])\in \mathrm{End}(V_{\rho}), \label{eqn:1-2}
 \end{multline}
and $r^0_{\rho}(u,v):=\langle u,v \rangle$. Then we have 
\begin{equation}
\sum_{0\le k \le N} m(\lambda_k)^{q}(p^{\rho}_{\lambda_k}(u))^{\ast}p^{\rho}_{\lambda_k}(v)=r^q_{\rho}(u,v),   \label{eqn:1-3}
\end{equation}
where $m(\lambda_k)$ is the conformal weight assigned from $V_{\rho}$ to $V_{\lambda_k}$. 
\end{thm}
In this paper, we will use the case of $q=0$ and $q=1$: 
\begin{gather}
\sum_{0\le k \le N}(p^{\rho}_{\lambda_k}(e_j))^{\ast}p^{\rho}_{\lambda_k}(e_i)=\delta_{ij},     \label{eqn-1} \\
\sum_{0\le k \le N} m(\lambda_k)(p^{\rho}_{\lambda_k}(e_j))^{\ast}p^{\rho}_{\lambda_k}(e_i)=-\frac{1}{4}\pi_{\rho}([e_j,e_i]). \label{eqn-2}
\end{gather}
\begin{rem}
The endomorphisms $\{ r^q_{\rho}(e_i,e_j)\}_{i,j}$ are useful to compute the eigenvalues of the higher Casimir operators (see \cite{NR} and \cite{Z}). 
\end{rem}
The Clifford homomorphisms also satisfy the following properties. 
\begin{prop}[\cite{H4}] \label{prop:1-2}
Let $u$ be in $\mathbf{R}^n$, $g$ in $Spin(n)$, and $[e_i,e_j]$ in $\mathfrak{spin}(n)$. Then we have 
\begin{equation}
p^{\rho}_{\lambda_k}(g u g^{-1})=\pi_{\lambda_k}(g)p^{\rho}_{\lambda_k}
(u)\pi_{\rho}(g^{-1}),         \label{eqn:2-17}
\end{equation}
and 
\begin{equation}
p_{\lambda_k}^{\rho}([[e_i,e_j],u])=\pi_{\lambda_k}([e_i,e_j])p_{\lambda_k}^{\rho}(u)-p_{\lambda_k}^{\rho}(u)\pi_{\rho}([e_i,e_j]). \label{eqn:2-18}
\end{equation}
\end{prop}
\section{Invariant operators on polynomials for higher bundles}\label{sec:2}
In the first part of this section, we give a well-known method to decompose the space of complex-valued polynomials on $\mathbf{R}^n$. We denote the canonical coordinate on $\mathbf{R}^n$ by $(x_i,\cdots,x_n)$, and the space of complex-valued polynomials with degree $q$ on $\mathbf{R}^n$ by $S^q$. The vector space $\sum_q S^q$ has the Hermitian inner product satisfying $(\partial /\partial x_i f(x),g(x))=(f(x),x_i g(x))$. The polynomial representation $(\pi_{s}, \sum S^q)$ of $\mathfrak{spin}(n)$ is defined by 
\begin{equation}
(\pi_{s}([e_k,e_l])f)(x):=4(-x_k\frac{\partial}{\partial x_l}+x_l\frac{\partial}{\partial x_k})f(x).
\end{equation}
To decompose the space $\sum S^q$, we use invariant operators compatible with the action of $\mathfrak{spin}(n)$. When the operator on $S^q$ maps to $S^{q-k}$, the order of the operator is said to be $k$. On $\sum S^q$, we have the following invariant operators: the Laplacian operator $\square:=-\sum \partial^2/\partial x_i^2$, and the $0$-th order operator $r\partial/\partial r=\sum x_i\partial/\partial x_i$ called the Euler operator, where $r^2$ is $\sum x_i^2$. The Euler operator measures the degree of polynomials. In other words, the vector space $S^q$ is the eigenspace with eigenvalue $q$ for the operator $r\partial/\partial r$. To decompose $S^q$ further, we use the $0$-th order invariant operator $-r^2\square$. This operator has the spectral decomposition corresponding to the irreducible decomposition. In fact, we show that $S^q$ is isomorphic to $\oplus_{0\le k \le [q/2]}r^{2k}H^{q-2k}$ and the eigenvalue of $-r^2\square$ on $r^{2k}H^{q-2k}$ is $k(2q-2k+n-2)$, where $H^q$ is the space of harmonic polynomials with degree $q$. In particular, the top component $H^q$ is the kernel of $-r^2\square$ and has the highest weight $h^q:=(q,0_{m-1})$. Thus, to decompose a representation space into irreducible components, we should investigate the spectral decompositions of invariant operators. 

Now, we shall consider the space of polynomials for higher bundles on $\mathbf{R}^n$. Let $(\pi_{\rho},V_{\rho})$ be an irreducible unitary representation of $\mathfrak{spin}(n)$. Then we have the (trivial) higher bundle $\mathbf{S}_{\rho}:=\mathbf{R}^n\times V_{\rho}$, and consider the polynomial sections of $\mathbf{S}_{\rho}$, that is, the $V_{\rho}$-valued polynomials $\sum S^q\otimes V_{\rho}$. This vector space is a representation space on where more invariant operators exist in addition to $-r^2\square$ and $r\partial/\partial r$. Here, the action of $\mathfrak{spin}(n)$ on $\sum_q S^q\otimes V_{\rho}$ is given as the tensor representation:
\begin{multline}
\mathfrak{spin}(n)\times S^q\otimes V_{\rho}\ni ([e_k,e_l], 
        f\otimes \phi )\to \\
  4(-x_k\frac{\partial}{\partial x_l}+x_l\frac{\partial}{\partial x_k})f \otimes\phi +  f\otimes \pi_{\rho}([e_k,e_l])\phi \in S^q\otimes V_{\rho}.
\end{multline}

We recall the Clifford homomorphism from $V_{\rho}$ to $V_{\lambda_k}$ given in Section \ref{sec:1}. By using the Clifford homomorphism, we introduce the following operators: 
\begin{gather}
x_{\lambda_k}^{\rho}:=\sum x_i p_{\lambda_k}^{\rho}(e_i) :S^q\otimes V_{\rho}\to S^{q+1}\otimes V_{\lambda_k}, \label{eqn:2-1} \\
(x_{\lambda_k}^{\rho})^{\ast}:=\sum x_i (p_{\lambda_k}^{\rho}(e_i))^{\ast}:S^q\otimes V_{\lambda_k}\to S^{q+1}\otimes V_{\rho}, \label{eqn:2-2} \\
D_{\lambda_k}^{\rho}:=\sum p_{\lambda_k}^{\rho}(e_i)\frac{\partial}{\partial x_i}:S^q\otimes V_{\rho}\to S^{q-1}\otimes V_{\lambda_k}, \label{eqn:2-3} \\
(D_{\lambda_k}^{\rho})^{\ast}:=-\sum (p_{\lambda_k}^{\rho}(e_i))^{\ast}\frac{\partial}{\partial x_i}:S^q\otimes V_{\lambda_k}\to S^{q-1}\otimes V_{\rho}. \label{eqn:2-4} 
\end{gather}
The differential operators $D_{\lambda_k}^{\rho}$ and $(D_{\lambda_k}^{\rho})^{\ast}$ are called the higher spin Dirac operators, which are generalization of the Dirac operator for higher bundles. If we define the inner product on $S^q\otimes V_{\rho}$ by the tensor inner product, then we show that the adjoint operators of $x_{\lambda_k}^{\rho}$ and $(x_{\lambda_k}^{\rho})^{\ast}$ are $-(D_{\lambda_k}^{\rho})^{\ast}$ and $D_{\lambda_k}^{\rho}$, respectively. 

We can show that the above operators are invariant operators on the $\mathfrak{spin}(n)$-module $\sum_q S^q\otimes V_{\rho}$.
\begin{prop}\label{prop:2-1}
The operators \eqref{eqn:2-1}-\eqref{eqn:2-4} are invariant operators. 
\end{prop}
\begin{proof}
We prove only the invariance of $x_{\lambda_k}^{\rho}$. It follows from the equation \eqref{eqn:2-18} that we have 
\begin{equation}
\begin{split}
 & ( -4x_k\frac{\partial}{\partial x_l}+4x_l\frac{\partial}{\partial x_k}+\pi_{\lambda_k}([e_k,e_l]) )x_{\lambda_k}^{\rho} \\ 
 =&\sum_i (-4x_k\frac{\partial}{\partial x_l}+4x_l\frac{\partial}{\partial x_k}+\pi_{\lambda_k}([e_k,e_l]) ) x_i p_{\lambda_k}^{\rho}(e_i) \\
=&\sum_i 4p_{\lambda_k}^{\rho}(e_i) (-\delta_{li}x_k-x_k x_i \frac{\partial}{\partial x_l} + \delta_{ki}x_l+ x_l x_i \frac{\partial}{\partial x_k} ) \\
   &\quad \quad \quad   +x_i \{p_{\lambda_k}^{\rho}(e_i)\pi_{\rho}([e_k,e_l]) +p_{\lambda_k}^{\rho}([[e_k,e_l],e_i])\} \\ 
=& 4(-p_{\lambda_k}^{\rho} (e_l)x_k+p_{\lambda_k}^{\rho}(e_k)x_l)+x_{\lambda_k}^{\rho}4(-x_k\frac{\partial}{\partial x_l}+x_l\frac{\partial}{\partial x_k})   \\
 &\quad \quad \quad +x_{\lambda_k}^{\rho}\pi_{\rho}([e_k,e_l])+\sum x_i(4\delta_{ki}p_{\lambda_k}^{\rho}(e_l)-4\delta_{li}p_{\lambda_k}^{\rho}(e_k)) \\
 =& x_{\lambda_k}^{\rho} ( -4x_k\frac{\partial}{\partial x_l}+4x_l\frac{\partial}{\partial x_k}+\pi_{\rho}([e_k,e_l]) ).
\end{split}
\end{equation}
\end{proof}

We shall investigate relations among these invariant operators, and reconstruct the Laplacian operator and the Euler operator. First, the formula \eqref{eqn-1} induces the following lemma. 
\begin{lem}\label{lem:2-2}
The invariant operators \eqref{eqn:2-1}-\eqref{eqn:2-4} satisfy that 
\begin{gather}
\sum_{0\le k \le N} (x_{\lambda_k}^{\rho})^{\ast}x_{\lambda_k}^{\rho}=\sum_i (x_i)^2=r^2, \quad
\sum_{0\le k \le N} (D_{\lambda_k}^{\rho})^{\ast}D_{\lambda_k}^{\rho}=\square, \\
\sum_{0\le k \le N} (D_{\lambda_k}^{\rho})^{\ast}x_{\lambda_k}^{\rho}=-n-r\frac{\partial}{\partial r}, \quad 
\sum_{0\le k \le N} (x_{\lambda_k}^{\rho})^{\ast}D_{\lambda_k}^{\rho}=r\frac{\partial}{\partial r}.
\end{gather}
\end{lem}
In similar way, the formula \eqref{eqn-2} gives the following lemma. 
\begin{lem}\label{lem:2-3}
The invariant operators \eqref{eqn:2-1}-\eqref{eqn:2-4} satisfy that 
\begin{equation}
\sum_{0\le k \le N} m(\lambda_k)(x_{\lambda_k}^{\rho})^{\ast}x_{\lambda_k}^{\rho}=0, \quad
\sum_{0\le k \le N} m(\lambda_k)(D_{\lambda_k}^{\rho})^{\ast}D_{\lambda_k}^{\rho}=0. \label{eqn-b}
\end{equation}
\end{lem}
\begin{rem}
The second equation in \eqref{eqn-b} means that $\mathbf{R}^n$ is a flat space (see \cite{H4}). 
\end{rem}
Since we have already given the decomposition of $S^q$, we shall decompose the $V_{\rho}$-valued harmonic polynomials $H^q\otimes V_{\rho}$. So we need relations among the Laplacian and the operators \eqref{eqn:2-1}-\eqref{eqn:2-4}. 
\begin{lem}\label{lem:2-4}
The Laplace operator $\square$ and the operators \eqref{eqn:2-1}-\eqref{eqn:2-4} satisfy that
\begin{gather}
[ \, \square, (D_{\lambda_k}^{\rho})^{\ast}]=0, \quad [\, \square, D_{\lambda_k}^{\rho}]=0, \\
[\, \square, x_{\lambda_k}^{\rho}]=-2D_{\lambda_k}^{\rho}, \quad [\, \square, (x_{\lambda_k}^{\rho})^{\ast}]=2(D_{\lambda_k}^{\rho})^{\ast}. 
\end{gather}
\end{lem}
From Lemma \ref{lem:2-3} and \ref{lem:2-4}, we have $0$-th order invariant operators compatible with the Laplacian $\square$. 
\begin{cor}\label{cor:2-5}
We consider the $0$-th order operators $\sum_{k}  m(\lambda_k)(D_{\lambda_k}^{\rho})^{\ast}x_{\lambda_k}^{\rho}$ and $\sum_{k} m(\lambda_k)(x_{\lambda_k}^{\rho})^{\ast}D_{\lambda_k}^{\rho}$. These operators commute with the Laplace operator:
\begin{equation}
[ \, \square, \sum_{k}  m(\lambda_k)(D_{\lambda_k}^{\rho})^{\ast}x_{\lambda_k}^{\rho}]=[\, \square, \sum_{k} m(\lambda_k)(x_{\lambda_k}^{\rho})^{\ast}D_{\lambda_k}^{\rho}]=0. 
\end{equation}
Furthermore, these two operators coincide with each other. 
\end{cor}
\begin{proof}
We can easily show that 
\begin{equation}
\begin{split}
 &\sum_k m(\lambda_k)(-(D_{\lambda_k}^{\rho})^{\ast}x_{\lambda_k}^{\rho}+(x_{\lambda_k}^{\rho})^{\ast}D_{\lambda_k}^{\rho})  \\ 
=&-\sum_{i,j} (x_j\frac{\partial}{\partial x_i}+x_i\frac{\partial}{\partial x_j})(\frac{1}{4}\pi_{\rho}([e_j,e_i])) \\
 =&0.
\end{split}\nonumber
\end{equation}
So we have proved the lemma.
\end{proof}
This corollary means that the operator $\sum_{k}  m(\lambda_k)(x_{\lambda_k}^{\rho})^{\ast}D_{\lambda_k}^{\rho}$ acts on $H^q\otimes V_{\rho}$ and has a spectral decomposition. 
\begin{prop}
Let $(\sum_{\mu} \pi_{\mu}, \sum_{\mu} V_{\mu})$ be the irreducible decomposition of $(\pi_{h^q} \otimes \pi_{\rho}, H^q \otimes V_{\rho})$. The $0$-th order invariant operator $\sum_{k} m(\lambda_k)(x_{\lambda_k}^{\rho})^{\ast} D_{\lambda_k}^{\rho}$ has the following spectral decomposition on $H^q\otimes V_{\rho}$: 
\begin{equation}
\sum_{k}  m(\lambda_k)(x_{\lambda_k}^{\rho})^{\ast}D_{\lambda_k}^{\rho}=m(\mu,q)   \quad \textrm{ on $V_{\mu}$}.
\end{equation}
The constant $m(\mu,q)$ is given by  
\begin{equation}
m(\mu,q):=\frac{1}{2}(q^2+(n-2)q+\|\rho+\delta\|^2-\|\mu+\delta \|^2), \label{eigen}
\end{equation}
where $\delta$ is half the sum of positive roots, and $\|\cdot\|$ is the canonical norm on the weight space, that is, $\| \nu \|^2=\sum_{1\le i \le m} (\nu^i)^2$. 
\end{prop}
\begin{proof}
We can show that
\begin{equation}
\begin{split}
 &\sum_k m(\lambda_k)(-(D_{\lambda_k}^{\rho})^{\ast}x_{\lambda_k}^{\rho}-(x_{\lambda_k}^{\rho})^{\ast}D_{\lambda_k}^{\rho}) \\
 =&-\sum_{ij} (-x_j\frac{\partial}{\partial x_i}+x_i\frac{\partial}{\partial x_j})(\frac{1}{4}\pi_{\rho}([e_j,e_i]))  \\
 =&-2\sum_{ij} \frac{1}{32}\pi_{h^q}([e_i,e_j])\otimes \pi_{\rho}([e_i,e_j]). 
 \end{split}
\end{equation}
The last equation is realized by using the Casimir operators. In fact, we can show that 
\begin{equation}
\sum_{ij} \frac{1}{32}\pi_{h^q}([e_i,e_j])\otimes \pi_{\rho}([e_i,e_j])
=C_{h^q\otimes \rho}-C_{h^q} \otimes \id-\id\otimes C_{\rho}.
\end{equation}
Here, the Casimir operator $C_{\nu}$ on the irreducible representation space $V_{\nu}$ is defined by 
\begin{equation}
C_{\nu}:=\frac{1}{64}\sum_{ij}\pi_{\nu}([e_i,e_j])\pi_{\nu}([e_i,e_j]),
\end{equation}
and acts as the constant $-(\|\delta+\nu\|^2-\|\delta\|^2)/2$ on $V_{\nu}$. 
Thus we have proved the proposition. 
\end{proof}
Instead of the $0$-th order operator in the above proposition, we consider the following operator corresponding to the Bochner type Laplacian on the bundle $\mathbf{S}_{\rho}$ (see \cite{H4}): 
\begin{equation}
  E:=\sum_{1 \le k \le N} (1-\frac{m(\lambda_k)}{m(\lambda_0)})(x^{\rho}_{\lambda_k})^{\ast}D^{\rho}_{\lambda_k},
\end{equation} 
where the weights $\{\lambda_k\}_k$ satisfy that $\lambda_0>\lambda_1>\cdots>\lambda_N$ with respect to the lexicographical order on the weight space. This operator $E$ is obtained by eliminateing the top operator $(x^{\rho}_{\lambda_0})^{\ast}D^{\rho}_{\lambda_0}$ from the equations $\sum_{k} (x_{\lambda_k}^{\rho})^{\ast}D_{\lambda_k}^{\rho}$ and $\sum_{k}  m(\lambda_k)(x_{\lambda_k}^{\rho})^{\ast}D_{\lambda_k}^{\rho}$. Then we have the following theorem. 
\begin{thm}\label{thm:2-7}
Let $(\sum_\mu \pi_{\mu}, \sum_\mu V_{\mu})$ be the irreducible decomposition of $(\pi_{h^q}\otimes \pi_{\rho}, H^q\otimes V_{\rho})$, where $h^q=(q,0_{m-1})$ and $\rho=(\rho^1,\cdots,\rho^m)$. The $0$-th order invariant operator $E$ is a non-negative operator and has the spectral decomposition on $H^q\otimes V_{\rho}$ as follows:
\begin{equation}
E=q+\frac{m(\mu,q)}{\rho^1} \quad \textrm{ on $V_{\mu}$}, 
\end{equation}
where the constant $m(\mu,q)$ is given in \eqref{eigen}. In particular, the $0$-eigenspace is the irreducible representation space with highest weight $\mu_0:=h^q+\rho$.
\end{thm}
In this theorem, we remark that the eigenvalues $\{e(\mu)\}$ of $E$ order as $0=e(\mu_0)<e(\mu_1)\le e(\mu_2)\le \cdots $ for $\mu_0>\mu_1\ge \mu_2\ge\cdots$. Here, the top component $(\pi_{\mu_0},V_{\mu_0})$ certainly exists with multiplication one. 
\begin{cor}
The irreducible representation with highest weight $\mu_0$ in $H^q \otimes V_{\rho}$ is realized as follows:
\begin{equation}
V_{\mu_0}=\bigcap_{1\le k \le N}\ker D^{\rho}_{\lambda_k}, 
\end{equation}
where $\ker D^{\rho}_{\lambda_k}$ is the kernel of $D^{\rho}_{\lambda_k}$ on $H^q \otimes V_{\rho}$.
 \end{cor}
\section{Examples}\label{sec:3}
In this section, we give some examples: spinor-valued harmonic polynomials and $p$-form-valued harmonic polynomials (see \cite{DSS}, \cite{Fo}-\cite{H2}, \cite{IT}, and \cite{T}).
\begin{exam}[spinor-valued harmonic polynomials]
We shall investigate only the odd dimensional case, that is, the case of $n=2m+1$. Let $V_{\Delta}$ be the spinor space with highest weight $\Delta=((1/2)_m)$. We consider the spinor-valued harmonic polynomials $H^q\otimes V_{\Delta}$, and have invariant operators: the Clifford multiplication $x=-x^{\ast}=\sum x_i e_i$ and the Dirac operator $D=D^{\ast}=\sum e_i \partial /\partial x_i$, twistor operator $T$ and so on. Then the $0$-th order invariant operator $E$ in Theorem \ref{thm:2-7} is $-x D=x^{\ast}D$. 

Now, we show that $H^q\otimes V_{\Delta}$ has the irreducible decomposition $V_{\mu_0}\oplus V_{\mu_1}$, where $\mu_0=h^q+\Delta=(q+1/2,(1/2)_{m-1})$ and $\mu_1=(q-1/2,(1/2)_{m-1})$. Then we have the spectral decomposition of $-x D$:
\begin{equation}
-x D=\begin{cases}
                0 &  \textrm{on $V_{\mu_0}$} \\ 
               n+2q-2  & \textrm{on $V_{\mu_1}$}.  \\
     \end{cases}
\end{equation}
In particular, we have 
\begin{equation}
V_{\mu_0}=\ker D, \qquad V_{\mu_1}=H^q\otimes V_{\Delta}/\ker D.
\end{equation}
\end{exam}
\begin{exam}[$p$-form-valued harmonic polynomials]
Let $\Lambda^p$ be the exterior tensor product space of $\mathbf{R}^n$ with degree $p$, which is the irreducible representation space with highest weight $(1_p,0_{m-p})$. We consider the $p$-form-valued harmonic polynomials $H^q\otimes \Lambda^p$, and have invariant operators: the exterior derivative $d=\sum e_i{}_{\wedge} \partial/\partial x_i$, its adjoint $d^{\ast}=-\sum i(e_i)\partial/\partial x_i$, the conformal killing operator $C$, $x_{\wedge}=\sum x_i e_i{}_{\wedge}$, and $i(x)=\sum x_i i(e_i)$ and so on. Here, $i(e_i)$ denotes the interior product of $e_i$. Then we have the spectral decomposition of $E=i(x)d-x_{\wedge}d^{\ast}$ on $H^q\otimes \Lambda^p$: 
\begin{equation}
i(x)d-x_{\wedge}d^{\ast}=\begin{cases}
                0 &  \textrm{on $V_{\mu_0}$} \\ 
                q+p & \textrm{on $V_{\mu_1}$}  \\
                n+q-p & \textrm{on $V_{\mu_2}$} \\
                n+2q-2 & \textrm{on $V_{\mu_3}$ (for $q\ge 2$)}, 
     \end{cases}
\end{equation}
where $\mu_0=(q+1,1_{p-1},0_{m-p})$, $\mu_1=(q,1_p,0_{m-p-1})$, $\mu_2=(q,1_{p-2},0_{m-p+1}$, and $\mu_3=(q-1,1_{p-1},0_{m-p})$. 
In particular, we have $V_{\mu_0}=\ker d\cap\ker d^{\ast}$.
\end{exam}
\section{Discussion}
In the case of $p$-form-valued harmonic polynomials, we can show that 
\begin{gather}
V_{\mu_1}=\ker d/\ker d\cap \ker d^{\ast}, \\
V_{\mu_2}=\ker d^{\ast}/\ker d\cap \ker d^{\ast}, \\
V_{\mu_3}=H^q\otimes \Lambda^p/(\ker d+\ker d^{\ast}).
\end{gather}
Thus, we can realize the irreducible components by using kernels of $d$ and $d^{\ast}$. In general case, we may realize any irreducible component of $H^q\otimes V_{\rho}$ by using kernels of higher spin Dirac operators (for the case of Rarita-Schwinger operator, see \cite{Bu}). 
\section*{Acknowledgements}
The author is partially supported by Waseda University Grant for Special Research Project 2000A-880.

\end{document}